\documentclass[a4paper,12pt]{amsart}
\usepackage{amsfonts}
\usepackage{amssymb}
\usepackage{ifthen}
\usepackage{graphicx}
\nonstopmode \numberwithin{equation}{section}
\usepackage[usenames]{color}
\newtheorem{thm}{Theorem}[section]
\newtheorem{lem}{Lemma}[section]
\newtheorem{cor}{Corollary}[section]

\newtheorem{cl}{Claim}[section]
\newtheorem{ca}{Case}[section]
\newtheorem{sca}{Subcase}[section]
\newtheorem{scl}[section]{Subclaim}
\newtheorem{conj}[equation]{Conjecture}

\theoremstyle{definition}
\newtheorem{defn}{Definition}[section]
\newtheorem{op}[equation]{Open Problem}
\newtheorem{ques}[equation]{Question}
\newtheorem{rem}{Remark}[section]
\newtheorem{exam}[equation]{Example}

\newcounter {own}
\def\theown {\thesection       .\arabic{own}}

\newenvironment{pf}[1][]{%
 \vskip 3mm
 \noindent
 \ifthenelse{\equal{#1}{}}%
  {{\slshape Proof. }}%
  {{\slshape #1.} }%
 }%
{\qed\bigskip}

\newcounter{alphabet}
\newcounter{tmp}

\makeatletter
\renewcommand{\Ref}[1]{\@ifundefined{r@#1}{}{\setcounter{tmp}{\ref{#1}}\Alph{tmp}}}
\makeatother
\newcommand{\IR}{{\mathbb R}}

\newcommand{\id}{{\operatorname{id}}}

\newcommand{\diam}{{\operatorname{diam}}}

\newcommand{\dist}{{\operatorname{dist}}}

\def\be{\begin{equation}}
\def\ee{\end{equation}}

\newcommand{\ben}{\begin{enumerate}}
\newcommand{\een}{\end{enumerate}}

\newcommand{\blem}{\begin{lem}}
\newcommand{\elem}{\end{lem}}
\newcommand{\bthm}{\begin{thm}}
\newcommand{\ethm}{\end{thm}}
\newcommand{\bcor}{\begin{cor}}
\newcommand{\ecor}{\end{cor}}
\newcommand{\beg}{\begin{exam}}
\newcommand{\eeg}{\end{exam}}
\newcommand{\begs}{\begin{examples}}
\newcommand{\eegs}{\end{examples}}
\newcommand{\bdefe}{\begin{defn}}
\newcommand{\edefe}{\end{defn}}
\newcommand{\bprob}{\begin{prob}}
\newcommand{\eprob}{\end{prob}}
\newcommand{\bques}{\begin{ques}}
\newcommand{\eques}{\end{ques}}
\newcommand{\bei}{\begin{itemize}}
\newcommand{\eei}{\end{itemize}}
\newcommand{\bcon}{\begin{conj}}
\newcommand{\econ}{\end{conj}}
\newcommand{\bop}{\begin{op}}
\newcommand{\eop}{\end{op}}

\newcommand{\bca}{\begin{ca}}
\newcommand{\eca}{\end{ca}}
\newcommand{\bsca}{\begin{sca}}
\newcommand{\esca}{\end{sca}}

\newcommand{\bcl}{\begin{cl}}
\newcommand{\ecl}{\end{cl}}

\newcommand{\bscl}{\begin{scl}}
\newcommand{\escl}{\end{scl}}

\newcommand{\bcons}{\begin{conjs}}
\newcommand{\econs}{\end{conjs}}
\newcommand{\bprop}{\begin{propo}}
\newcommand{\eprop}{\end{propo}}
\newcommand{\br}{\begin{rem}}
\newcommand{\er}{\end{rem}}
\newcommand{\brs}{\begin{rems}}
\newcommand{\ers}{\end{rems}}
\newcommand{\bo}{\begin{obser}}
\newcommand{\eo}{\end{obser}}
\newcommand{\bos}{\begin{obsers}}
\newcommand{\eos}{\end{obsers}}
\newcommand{\bpf}{\begin{pf}}
\newcommand{\epf}{\end{pf}}
\newcommand{\ba}{\begin{array}}
\newcommand{\ea}{\end{array}}
\newcommand{\beq}{\begin{eqnarray}}
\newcommand{\beqq}{\begin{eqnarray*}}
\newcommand{\eeq}{\end{eqnarray}}
\newcommand{\eeqq}{\end{eqnarray*}}

\newcommand{\ds}{\displaystyle}

\newcounter{minutes}
\divide\time by 60
\newcounter{hours}
\multiply\time by 60 \addtocounter{minutes}{-\time}
\begin{document}

\bibliographystyle{amsplain}
\title{Teichm\"{u}ller's problem for Gromov hyperbolic domains}

\author{Qingshan Zhou${}^{\mathbf{*}}$}
\address{Qingshan Zhou, School of Mathematics and Big Data, Foshan university,  Foshan, Guangdong 528000, People's Republic
of China} \email{qszhou1989@163.com; q476308142@qq.com}

\author{Antti Rasila}
\address{Antti Rasila, Technion -- Israel Institute of Technology, Guangdong Technion, Shantou, Guangdong 515063, People's Republic
of China} \email{antti.rasila@gtiit.edu.cn; antti.rasila@iki.fi}

%%%%%%%% BEGIN TIMESTAMP
\def\thefootnote{}
\footnotetext{ \texttt{\tiny File:~\jobname .tex,
          printed: \number\year-\number\month-\number\day,
          \thehours.\ifnum\theminutes<10{0}\fi\theminutes}
} \makeatletter\def\thefootnote{\@arabic\c@footnote}\makeatother
%%%%%%%% END TIMESTAMP

\date{}
 \keywords{
Teichm\"{u}ller's problem, quasiconformal, Gromov hyperbolic domains, quasi-hyperbolic metric, distance ratio metric.\\
${}^{\mathbf{*}}$ Corresponding author}

\begin{abstract} Let $\mathcal{T}_K(D)$ be the class of $K$-quasiconformal automorphisms of a domain $D\subsetneq \mathbb{R}^n$ with identity boundary values. Teichm\"uller's problem is to determine how far a given point $x\in D$ can be mapped under a mapping $f\in \mathcal{T}_K(D)$. We estimate this distance between $x$ and $f(x)$ from the above by using two different metrics, the distance ratio metric and the quasihyperbolic metric. We study Teichm\"{u}ller's problem for Gromov hyperbolic domains in $\mathbb{R}^n$ with identity values at the boundary of infinity. As applications, we obtain results on Teichm\"{u}ller's problem for $\psi$-uniform domains and inner uniform domains in $\mathbb{R}^n$.
\end{abstract}

\thanks{
Qingshan Zhou was supported by NSF of China (No. 11901090), and by Department of Education of Guangdong Province, China (Grant Nos. 2018KQNCX285 and 2018KTSCX245). Antti Rasila was supported by NSF of China (No. 11971124).}

\maketitle{} \pagestyle{myheadings} \markboth{}{}

\section{Introduction and main results}
Teichm\"{u}ller's problem concerns finding a lower bound for the maximal dilation of the class of quasiconformal self-maps of a domain $D$, with identity boundary values, moving a point $x$ in the domain to a given point. Following \cite{VZ}, suppose that $D\subsetneq \mathbb{R}^n$ ($n\geq 2$) is a domain. Note that the boundary of a domain in $\mathbb{R}^n$ is taken in the topology of the Riemann sphere $\widehat{\mathbb{R}^n}=\mathbb{R}^n\cup \{\infty\}$, so $\partial D$ contains at least two points. Let
$$\mathcal{T}_K(D)=\Big\{\overline{D} \xrightarrow{f}  \overline{D}\;\big  |\; f\; \mbox{is a homeomorphism s.t.} \;f|_D\;  \mbox{is K-QC and}\;f|_{\partial D}=\id_{\partial D}\Big\},$$
where $K$-QC means $K$-quasiconformal and $\id_{\partial D}$ denotes the identity map on $\partial D$.

In \cite{T}, Teichm\"{u}ller considered the above class of maps $\mathcal{T}_K(D)$ with $D=\mathbb{R}^2\setminus\{(0,0),(1,0)\}$, and obtained the following sharp inequality:
$$h_D(x,f(x))\leq \log K$$
for all $x\in D$, where $h_D$ is the hyperbolic metric of $D$. This result may be regarded as a stability result for quasiconformal homeomorphisms, which hold the boundary pointwise fixed and map the domain onto itself.

For Teichm\"{u}ller type results concerning the same problem in the case of the unit balls of $\mathbb{R}^n$, $n\ge 2$, we refer to \cite{AV78,KTV,K68,LVW,MV11,M09}. Vuorinen and Zhang have further studied Teichm\"{u}ller's problem for other domains in $\mathbb{R}^n$, such as convex domains and uniform domains with uniformly perfect boundaries, see \cite{Vu84, VZ}. We also note that Teichm\"{u}ller type results  are applicable to questions related to the homogeneity of domains, see \cite{GP,KM11} and references therein.

In this paper, we investigate Teichm\"{u}ller's problem for domains in $\mathbb{R}^n$ with uniformly perfect (see Definition \ref{unif-perf}) boundaries. We prove that the distance between $x$ and its quasiconformal image $f(x)$ is uniformly bounded with respect to the distance ratio metric $j$. For the definition of the metric $\mathit{j}$ see (\ref{ss-1}).

\begin{thm}\label{thm-1}
Let $n\geq 2, C>1$, and $K\ge1$. There exists a constant $M=M(n,C,K)$ such that:
If
\begin{enumerate}
  \item $D\subsetneq \mathbb{R}^n$ is a domain,
  \item $\partial D$ is a $C$-uniformly perfect set,
  \item $f \in \mathcal{T}_K(D)$,
\end{enumerate}
then $j_D(x,f(x))\leq M$ for all $x\in D$.

\end{thm}

\br\label{a-2} Let $D \subsetneq \mathbb{R}^n$ be a domain. We say that $\partial D$ is $C$-uniformly perfect, if $\partial D$ is  $C$-uniformly perfect with respect to the spherical metric. It follows from \cite[Lemma C]{WZ17} that uniform perfectness is preserved by quasim\"obius transformations, and thus $\partial D$ is a $C'$-uniformly perfect set with respect to the Euclidean metric of $\mathbb{R}^n$, where $C'$ depends only on $C$.

Note that Theorem \ref{thm-1} does not hold for domains in $\mathbb{R}^n$ whose boundary is not uniformly perfect.  Let $o=(0,0)$, $D=B(o,1)\setminus \{o\}\subset \mathbb{R}^2$, and $x_m=(1/m,0)$ for $m\geq 3$. Consider the function $f:\mathbb{R}^2\to \mathbb{R}^2$ defined by $f(x)=|x|x$ for all $x\in \mathbb{R}^2$. Therefore it follows from \cite[5.21]{Va99} that $f\in \mathcal{T}_K(D)$ with $K=2$. However
$$j_D(x_m,f(x_m))\to\infty,$$
as $m\to \infty$.
\er

\br$($\cite[Remark 3.5]{Vu84}$)$ Next we give an example to show that the constant $M$ of Theorem \ref{thm-1} is strictly greater than $\log 3$ even if $f$ is a conformal map. Let $D=\mathbb{R}^3\setminus Z$ and $Z=\{(0,0,z):z\in \mathbb{R}\}$. Let $e_1=(1,0,0)$, and let $f$ be the rotation around the line $Z$ with $f(e_1)=-e_1$ such that $f$ keeps the  line $Z$ pointwise fixed. Then $f$ is conformal. However, we have $j_D(x,f(x))=\log 3$ for all $x\in D$.
\er

We may apply Theorem \ref{thm-1} to  Teichm\"{u}ller's problem for $\psi$-uniform domains (see Definition \ref{i-3}) and deduce that the quasihyperbolic distance $k$ (see Definition \ref{i-2}) between $x$ and its image point $f(x)$ is uniformly bounded.

\bcor\label{cor-1}  Let $n\geq 2, C>1$, $K\ge1$, and let $\psi:[0,\infty)\to [0,\infty)$ be a homeomorphism. There exists a constant $M'=M'(n,C,K,\psi)$ such that:
If
\begin{enumerate}
  \item $D\subsetneq \mathbb{R}^n$ is a $\psi$-uniform domain,
  \item $\partial D$ is a $C$-uniformly perfect set,
  \item $f \in \mathcal{T}_K(D)$,
\end{enumerate}
then $k_D(x,f(x))\leq M'$ for all $x\in D$.

\ecor

\br The class of $\psi$-uniform domains in $\mathbb{R}^n$ was introduced by Vuorinen in \cite{Vu85}. In \cite{HKSV}, H\"{a}st\"{o}, Kl\'{e}n, Sahoo and Vuorinen studied certain geometric properties of these domains. Because convex domains and uniform domains are both $\psi$-uniform,  Corollary \ref{cor-1} holds for all convex domains and for all uniform domains in $\mathbb{R}^n$.
\er
Recently, Bonfert-Taylor, Canary, Martin and Taylor studied Teichm\"{u}ller's problem in the case of the classical hyperbolic space $\mathbb{H}^n$, and they proved that if the boundary extension of a quasiconformal map is the identity on $\partial \mathbb{H}^n$, then it is uniformly close to the identity map on $\mathbb{H}^n$; see \cite[Lemma 4.1]{BCMT05}. In \cite{BHK}, Bonk, Heinonen and Koskela proved that every uniform domain  in $\mathbb{R}^n$ is Gromov hyperbolic with respect to the quasihyperbolic metric. This motivates us to study Teichm\"{u}ller's problem for Gromov hyperbolic domains.

As the second main aim of this paper, we consider Teichm\"{u}ller's problem for Gromov hyperbolic domains in $\mathbb{R}^n$. Let $D\subsetneq \mathbb{R}^n$ ($n\geq 2$) be a domain and $k_D$ its quasihyperbolic metric. If $(D,k_D)$ is a $\delta$-hyperbolic metric space for some $\delta\geq 0$, then we call $D$ a {\it Gromov hyperbolic domain} or a $\delta$-{\it hyperbolic domain}. Denote by $\partial^*D$ the Gromov boundary of the hyperbolic space $(D,k_D)$ and by $D^*$ its Gromov closure. For more information about Gromov hyperbolic spaces see Subsection \ref{sec-g}.

Let $D\subsetneq \mathbb{R}^n$ be a $\delta$-hyperbolic domain and let $f:D\to D$ be a $K$-quasiconformal homeomorphism. Note that $(D,k_D)$ is a proper geodesic metric space by \cite[Proposition 2.8]{BHK}, and it is not difficult to see from \cite[Theorem 3]{GO} that $f:(D,k_D)\to (D,k_D)$ is a rough quasi-isometry. By combining these two facts with \cite[Proposition 6.3]{BS}, we find that the image of any Gromov sequence under $f$ is also Gromov, and so, $f$ induces a boundary map $\partial f:\partial^*D\to \partial^*D$. Now define
$$\mathcal{T}_{K}^*(D) =\Big\{D \xrightarrow{f} D\;\big|\;f\;\mbox{is}\; K\mbox{-QC so that}\;\partial f=\id_{\partial^* D} \Big\},$$
where $f|_D$ is $K$-quasiconformal with respect to the Euclidean metric of $\mathbb{R}^n$. Our second main result reads as follows:

\begin{thm}\label{thm-2} Let $n\geq 2$, $\delta\geq0$, $C>1$, and $K\ge1$. There exists a constant $L=L(n,\delta,C,K)$ such that:
If
\begin{enumerate}
  \item $D\subsetneq \mathbb{R}^n$ is a $\delta$-hyperbolic domain,
  \item $\partial^* D$ equipped with a visual metric $\rho$ is $C$-uniformly perfect,
  \item $f \in \mathcal{T}_K^*(D)$,
\end{enumerate}
then  $k_D(x,f(x))\leq L$ for all $x\in D$.
\end{thm}
\br The definition of visual metrics on $\partial^* D$ is given by  Definition \ref{a-1}. By \cite[Corollary 5.2.9]{BuSc}, we see that $\partial^* D$ endowed with any two visual metrics are quasim\"obius equivalent.
It follows from \cite[Lemma C]{WZ17} that uniform perfectness is a quasim\"obius invariant and thus $\partial^* D$ is $C'$-uniformly perfect with respect to any visual metric, where $C'$ depends only on $C$, $\delta$ and the parameters of the visual metrics.

\er

It follows from \cite[Theorem 6.5]{BS} that if two Gromov hyperbolic spaces are roughly isometrically equivalent, then their boundaries at infinity are bilipschitzly equivalent with respect to visual metrics based on the same parameters. Conversely, if the boundaries at infinity of two roughly starlike Gromov hyperbolic geodesic spaces are bilipschitz equivalent, then these spaces are roughly isometrically equivalent; see \cite[Theorem 7.1.2]{BuSc}. Indeed, we know from Theorem \ref{thm-2} that a quasiconformal map $f\in \mathcal{T}_{K}^*(D)$ is a rough isometry with respect to the quasihyperbolic metric.

\bcor\label{cor-2} Let $n\geq 2$, $\delta\geq0$, $C>1$, and $K\ge1$. There exists a constant $L'=L'(n,\delta,C,K)$ such that:
If
\begin{enumerate}
  \item $D\subsetneq \mathbb{R}^n$ is a $\delta$-hyperbolic domain,
  \item $\partial^* D$ is a $C$-uniformly perfect set,
  \item $f \in \mathcal{T}_K^*(D)$,
\end{enumerate}
then  for all $x,y\in D$,
$$k_D(x,y)-L'\leq k_D(f(x),f(y)) \leq k_D(x,y) +L'.$$
\ecor

As the second application of Theorem \ref{thm-2}, we investigate Teichm\"{u}ller's problem for inner uniform domains $D$ in $\mathbb{R}^n$. For $x,y\in D$, the {\it inner Euclidean metric $d_I$} of $D$ is given by
$$d_I(x,y):=\inf\{\ell(\gamma_{x,y})\},$$
where the infimum is taken over all rectifiable curves $\gamma_{x,y}$ in $D$ with endpoints $x$ and $y$.

Let $D_I=(D,d_I)$ and let $\overline{D}_I$ be the metric completion of $D$ with respect to $d_I$. Following \cite[Section 6]{HSX}, let $\widehat{D}_I$ be the one point compactification $\overline{D}_I\cup\{\infty\}$ of $\overline{D_I}$ if $\overline{D}_I$ is unbounded, and $\widehat{D}_I=\overline{D}_I$ if $\overline{D}_I$ is bounded. Denote by $\partial_I D$ the topological boundary of $D_I$ in $\widehat{D}_I$. Thus $\partial_I D=\overline{D}_I\setminus D$ if $(D,d_I)$ is bounded, and $\partial_I D=(\overline{D}_I\cup\{\infty\})\setminus D$ if $(D,d_I)$ is unbounded. Now define
$$\mathcal{T}_{K}(D_I)=\Big\{ \widehat{D}_I \xrightarrow{f} \widehat{D}_I \;\big|\; f\mbox{ is a homeomorphism s.t.}\; f|_{D} \;\mbox{is K-QC and}\;f|_{\partial_I D}=\id_{\partial_I D} \Big\},$$
where $f|_D$ is $K$-quasiconformal with respect to the Euclidean metric of $\mathbb{R}^n$. Our result concerning Teichm\"{u}ller's problem for inner uniform domains (see Definition \ref{i-1}) is the following:

\bcor\label{cor-3} Let $n\geq 2, A\geq1, C>1$, and $K\ge1$. There exists a constant $H=H(n,A,C,K)$ such that:
If
\begin{enumerate}
  \item $D\subsetneq \mathbb{R}^n$ is an inner $A$-uniform domain,
  \item $\partial_I D$ is a $C$-uniformly perfect set,
  \item $f\in \mathcal{T}_{K}(D_I)$,
\end{enumerate}
then  for all $x\in D$,
$$k_D(x,f(x))\leq H.$$
\ecor

\br
We are grateful to Matti Vuorinen for pointing out that the main motivation in studying Teichm\"{u}ller's problem is to find sharp estimates, as was already done in Teichm\"{u}ller's original work \cite{T}. More precisely, it is interesting to understand the convergence behavior of the bounds in the stability theory whenever the quasiconformality coefficient $K$ tends to $1$, see, e.g., \cite{KTV,MV11,VZ}. In light of the proofs of Theorems \ref{thm-1} and \ref{thm-2}, the related bounds are unlikely to be sharp.

Indeed, in studying problems on Gromov hyperbolic spaces, properties of geometries at large scales are the main area of concern. In particular, in order to show that any quasiconformal self-map of a uniformly quasiconformally homogeneous manifold is uniformly close to an isometry, Bonfert-Taylor, Canary, Martin, and Taylor \cite{BCMT05} only needed a suitable bound, but sharpness was not required in their study of Teichm\"uller's problem for hyperbolic manifolds.
\er

The rest of this paper is organized as follows.  In Section \ref{sec-2}, we recall necessary definitions and preliminary results. The proof of Theorem \ref{thm-1} is given in Section \ref{sec-3}. Section \ref{sec-4} is devoted to the proof of Theorem \ref{thm-2}, and the proofs of Corollaries \ref{cor-1}, \ref{cor-2} and \ref{cor-3} are presented in Section \ref{sec-5}.

\section{Preliminaries and auxiliary results}\label{sec-2}
\subsection{Notation}
Let letter $A,B,C,\ldots$ denote positive numerical constants. Similarly, $C(a,b,c,\ldots)$ denotes universal positive functions of the parameters $a,b,c,\ldots$. Sometimes we write $C=C(a,b,c,\ldots)$ to emphasize the parameters on which $C$ depends and abbreviate $C(a,b,c,\ldots)$ to $C$.
\subsection{Metric geometry}\label{subsec-2.1}

Let $(X, |\cdot|)$ be a metric space, and let $$B(x,r)=\big\{ z\in X \;\big|\; |z-x|< r\big\}.$$

The metric space $X$ is called {\it proper} if its closed balls are compact. For a bounded set $S\subseteq X$, we denote the diameter of $S$ by $\diam (S)$. We use $\overline{X}$ to denote the metric completion of $X$ and $\partial X=\overline{X}\setminus X$ to be its metric boundary. A metric space $X$ is called {\it incomplete} if it is not complete. Thus incompleteness of $X$  implies that $\partial X\neq \emptyset$. The identity map of $X$ is denoted by $\id_X$.

A domain $D\subseteq X$ is an open and connected non-empty set. Let $X$ be a connected and complete metric space, and let $D \subsetneq X$ be a domain. Thus $\partial D\neq \emptyset$, and we write
$$d(x)=\dist(x,\partial D)$$
for all $x\in D$. For $x,y\in D$, the {\it distance ratio distance} $j_D(x,y)$ is defined by
\be\label{ss-1} j_D(x,y)=\log\Big(1+\frac{|x-y|}{\min\{d(x),d(y)\}}\Big).\ee

\bdefe
\label{unif-perf}
Let $C>1$. A metric space $X$ is {\it $C$-uniformly perfect}, if for each $x\in X$ and every $r>0$, $B(x,r)\setminus B(x, r/C)\not=\emptyset$ provided $X\setminus B(x,r)\not=\emptyset$.
\edefe

We also record the following invariance property of uniform perfectness of metric spaces under quasim\"obius maps (see Definition \ref{d-qm}) for later use.

\begin{lem}\label{Thm-z0} $($\cite[Lemma C]{WZ17}$)$ Suppose that $f:X\to Y$ is a $\theta$-quasim\"obius homeomorphism between two metric spaces. If $X$ is $C$-uniformly perfect, then $Y$ is $C'$-uniformly perfect with $C'=C'(C,\theta)$.
\end{lem}

A curve is a continuous function $\gamma:$ $\mathbb{R}\supset[a,b]\to X$. The length of $\gamma$ is defined by
$$\ell(\gamma)=\sup\Big\{\sum_{i=1}^{n}|\gamma(t_i)-\gamma(t_{i-1})|\Big\},$$
where the supremum is taken over all partitions $a=t_0<t_1<t_2<\ldots<t_n=b$. The curve $\gamma$ is called {\it rectifiable} if $\ell(\gamma)<\infty$. The metric space $X$ is called {\it rectifiably connected} if each pair of points can be connected by a rectifiable curve.

The length function associated with a rectifiable curve $\gamma$: $[a,b]\to X$ is $s_{\gamma}$: $[a,b]\to [0, \ell(\gamma)]$, defined by
$s_{\gamma}(t)=\ell(\gamma|_{[a,t]})$ for $t\in [a,b]$. For any rectifiable curve $\gamma:$ $[a,b]\to X$, there is a unique parametrization $\gamma_s:$ $[0, \ell(\gamma)]\to X$ such that $\gamma=\gamma_s\circ s_{\gamma}$. Obviously, $\ell(\gamma_s|_{[0,t]})=t$ for $t\in [0, \ell(\gamma)]$. The parametrization $\gamma_s$ is called the {\it arclength parametrization} of $\gamma$. For a rectifiable curve $\gamma$ in $X$, following \cite[Section 10]{BHK}, the line integral over $\gamma$ of each Borel function $\varrho:$ $X\to [0, \infty)$ is
$$\int_{\gamma}\varrho \;ds=\int_{0}^{\ell(\gamma)}\varrho\circ \gamma_s(t) \;dt.$$

In 1978, uniform domains in $\mathbb{R}^n$ were introduced by Martio and Sarvas \cite{MS}. In order to establish their uniformization theory of Gromov hyperbolic spaces, Bonk, Heinonen and Koskela \cite{BHK} generalized this concept to the setting of metric spaces.

\bdefe
Let $A\geq 1$, and let $X$ be a rectifiably connected and complete metric space. A domain $D\subsetneq X$ is called $A$-{\it uniform} if each pair of points $x$, $y$ in $D$ can
be joined by a rectifiable arc $\gamma$ in $D$ satisfying:
\begin{enumerate}
\item $\ell(\gamma)\leq A\,|x-y|$, and
\item $\min\{\ell(\gamma[x,z]),\ell(\gamma[z,y])\}\leq A\,d(z)$ for all $z\in \gamma$,
\end{enumerate}
\noindent where $\gamma[x,z]$ is the part of $\gamma$ between $x$ and $z$.
\edefe

We remark that $D$ is called $A$-{\it quasiconvex}, if any two points of $D$ can be connected by a curve satisfying the condition $(1)$ above. We call a domain {\it uniform} if it is $A$-uniform for some constant $A\geq 1$ and {\it quasiconvex} if it is $A$-quasiconvex for some $A\geq 1$.

\bdefe\label{i-1}
We say that a domain $D \subsetneq\mathbb{R}^n$ is {\it inner uniform}, if $(D,d_I)$ is $A$-uniform for some $A\geq 1$, where $d_I$ is the inner Euclidean metric $d_I$ of $D$.
\edefe

\bdefe\label{i-2}
Let $X$ be a connected and complete metric space, and let $D\subsetneq X$ be a rectifiably connected domain. The {\it quasihyperbolic metric} $k_D$ of $D$ is defined as
$$k_D(x,y)=\inf \int_\gamma \frac{|dz|}{d(z)},$$
where the infimum is taken over all rectifiable curves $\gamma$ joining $x$ and $y$ in $D$.
\edefe

There is an important property of uniform domains associated to $k_D$ and the distance ratio metric $j_D$. The statement is as follows.

\begin{lem}\label{Thm-z1} $($\cite[Lemma 2.13]{BHK}$)$
Let $X$ be a locally compact, rectifiably connected metric space and $D\subsetneq X$ an $A$-uniform domain. Then for all $x,y\in D$,
$$ k_D(x,y)\leq 4A^2j_D(x,y).$$
\end{lem}

\bdefe\label{i-3}
Let $\psi:[0,\infty)\to [0,\infty)$ be a homeomorphism. A domain $D\subsetneq \mathbb{R}^n$ is called $\psi$-{\it uniform} if for all $x$, $y$ in $D$,
$$k_D(x,y)\leq \psi(r_D(x,y))\;\;\;\;\;\mbox{where}\;\;\;\;\;r_D(x,y)=\frac{|x-y|}{\min\{d(x),d(y)\}}.$$
\edefe

Note that an $A$-uniform domain is $\psi$-uniform with $\psi(t)=4A^2\log(1+t)$, which follows from Lemma \ref{Thm-z1}.

\subsection{Maps on metric spaces} Assume that $X$ and $Y$ are metric spaces.
For the basic theory of quasiconformal maps we refer to \cite{Hei,Vaibook}. There are several equivalent definitions for quasiconformality in $\mathbb{R}^n$. We adopt a version of the metric definition.

\bdefe
Let $n\geq 2$, let $D$ and $D'$ be domains in $\mathbb{R}^n$, and let $f:D\to D'$ be a homeomorphism. For $1\leq K< \infty$, we say that $f$ is $K$-{\it quasiconformal} if
$$H(x):=\limsup_{r\rightarrow 0}\frac{\sup\big\{|f(x)-f(y)|\; |\;|x-y|= r\big\}}{\inf\big\{|f(x)-f(z)|\;  |\;|x-z|= r\big\}}\leq K$$
for all $x\in D$, and that $f$ is {\it quasiconformal} if it is $K$-quasiconformal for some $K$.
\edefe

For $K$-quasiconformal maps we have the following property:

\begin{thm}\label{Thm-1} $($\cite[Theorem 3]{GO}$)$
For $n\geq 2$, $K\geq 1$, there are constants $C\geq 1,\mu\in(0,1]$ depending only on $n$
and $K$ such that if $D, D'\subsetneq \mathbb{R}^n$ and
$f: D\to D'$ is a $K$-quasiconformal map, then for all $x,y\in D$,
$$k_{D'}(f(x),f(y))\leq C \max\{k_D(x,y), k_D(x,y)^{\mu}\}.$$
\end{thm}

Following notations and terminology of \cite{Hei,HL,TV,Va85,Va99,Vu88}, we next recall the definitions of quasisymmetric and quasim\"obius maps:

\bdefe A homeomorphism $f$ from $X$ to $Y$ is said to be
\begin{enumerate}
\item $\eta$-{\it quasisymmetric} if there is a homeomorphism $\eta : [0,\infty) \to [0,\infty)$ such that
$$ |x-a|\leq t|x-b|\;\; \mbox{implies}\;\;   |f(x)-f(a)| \leq \eta(t)|f(x)-f(b)|$$
for each $t>0$ and for each triplet $x,$ $a$, $b$ of points in $X$;

\item {\it weakly $H$-quasisymmetric} if there is a constant $H<\infty$ such that
$$ |x-a|\leq |x-b|\;\;  \mbox{ implies}\;\;   |f(x)-f(a)| \leq H|f(x)-f(b)|$$
for each triplet $x$, $a$, $b$ of points in $X$.

\end{enumerate}\edefe

\bdefe Let $X$ and $Y$ be incomplete and connected metric spaces. Let $0<q<1$, $H\geq 1$,  and $\eta : [0,\infty) \to [0,\infty)$ a homeomorphism. Suppose $f:X\to Y$ is a homeomorphism.

The map $f$ is said to be {\it $q$-locally $\eta$-quasisymmetric} if the restrictions $f|_{B(z,q d(z))}$ of $f$ to $B(z,q d(z))$ are $\eta$-quasisymmetric for all $z\in X$, where $d(z)=\dist(z,\partial X)$.

Similarly, $f$ is called {\it $q$-locally weakly $H$-quasisymmetric} if the restrictions $f|_{B(z,q d(z))}$ of $f$ to $B(z,q d(z))$ are weakly $H$-quasisymmetric for all $z\in X$.
\edefe

The following result is needed in the proof of Theorem \ref{thm-1}.

\begin{thm}\label{Thm-A} $($\cite[Theorem 3.12]{Va90}$)$ Let $n\geq 2$, $K\geq 1$ and $\eta : [0,\infty) \to [0,\infty)$ a homeomorphism. Suppose that $D$ and $D'$ are domains in $\mathbb{R}^n$ and $f:\overline{D} \to \overline{D'}$ is a homeomorphism such that $f|_D$ is $K$-quasiconformal and $f|_{\partial D}$ is $\eta$-quasisymmetric. Then $f$ is $\eta_1$-quasisymmetric with $\eta_1=\eta_1(K,\eta,n)$.
\end{thm}
%%%%%%%%%%%%%%%%%
A quadruple in $X$ is an ordered sequence $Q = (a, b, c, d)$ of four distinct points
in $X$. The {\it cross ratio} of $Q$ is defined to be the number
$$\tau(Q) = |a, b, c, d| = \frac{|a-c|}{|a-b|}\cdot
\frac{|b-d|}{|c-d|}.$$

\bdefe\label{d-qm} Note that a homeomorphism $f$ from $X$ to $Y$ is said to be {\it $\theta$-quasim\"obius} if $\theta : [0,\infty) \to [0,\infty)$ is a homeomorphism such that
$$\tau(f(Q))\leq\theta(\tau(Q))$$ holds for each quadruple $Q\subseteq X$.
\edefe

There is a criterion for quasim\"{o}bius maps between two bounded metric spaces to be quasisymmetric given by V\"{a}is\"{a}l\"{a}. For later use we record this result as follows.

\begin{thm}\label{Thm-2} $($\cite[Theorem 3.12]{Va85}$)$
Suppose that $X$ and $Y$ are two bounded metric spaces, that $\lambda>1$, that $z_1,$ $z_2,$ $z_3$ in $X$, and that $f: X\to Y$ is a $\theta$-quasim\"obius homeomorphism satisfying the {\it three-point condition}:
\be\label{a-6} |z_i-z_j|\geq \frac{1}{\lambda}\diam(X)\;\;\mbox{and}\;\; |f(z_i)-f(z_j)|\geq \frac{1}{\lambda}\diam(Y)\ee
for all $i\neq j\in\{1,2,3\}$. Then $f$ is $\eta$-quasisymmetric with $\eta=\eta(\theta,\lambda)$.
\end{thm}

\bdefe
Let $f: X\to Y$ be a map (not necessarily continuous) between metric spaces $X$ and $Y$, and let $L\geq 1$ and $M\geq 0$ be constants. \begin{enumerate}
\item
If \begin{enumerate}
     \item for each $x'\in Y$, there is $x\in X$ with $|x'-f(x)|\leq M$, and
     \item for all $x,y\in X$,
$$L^{-1}|x-y|-M\leq |f(x)-f(y)|\leq L|x-y|+M,$$
   \end{enumerate}

then $f$ is called an {\it $(L, M)$-roughly quasi-isometric map} (cf. \cite{BS}).
If $L=1$, then $f$ is called an {\it $M$-roughly isometric map}.
\noindent
\item Moreover, if $f$ is a homeomorphism and $M=0$, then it is called an {\it $L$-bilipschitz map}.
\end{enumerate}\edefe

\subsection{Gromov hyperbolic spaces}\label{sec-g}

In this subsection, we recall some necessary terminology concerning Gromov hyperbolic spaces (cf. \cite{BHK,BS,BrHa,BuSc}). Let $(X,d)$ be a metric space. Fix a base point $w$ in $X$.
\begin{enumerate}
\item
For $x,y\in X$, let
$$(x|y)_w=\frac{1}{2}\big(d(x,w)+d(y,w)-d(x,y)\big).$$
This number is called the  {\it Gromov product} of $x,y$ with respect to $w$.
\item The space $X$ is called {\it geodesic}, if each pair of points $x,y\in X$ can be joined by a geodesic $[x,y]$; that is, a curve whose length is precisely the distance between $x$ and $y$. Moreover, a {\it geodesic triangle} $\Delta$ is a set $\Delta=[x_1,x_2]\cup[x_2,x_3]\cup[x_3,x_1]\subseteq X$.

\item
Suppose $(X,d)$ is geodesic. The metric space $X$ is called {\it $\delta$-hyperbolic} $(\delta\geq 0)$ if each point on the edge of any geodesic triangle in $X$ is within distance $\delta$ of some point on one of the other two edges. If $X$ is $\delta$-hyperbolic for some $\delta\geq 0$, we also say that it is Gromov hyperbolic.
\item
Suppose $(X, d)$ is $\delta$-hyperbolic.
\begin{enumerate}
\item
A sequence $\{x_i\}$ in $X$ is called a {\it Gromov sequence} if $(x_i|x_j)_w\rightarrow \infty$ as $i,$ $j\rightarrow \infty.$
\item
Two such sequences $\{x_i\}$ and $\{y_j\}$ are said to be {\it equivalent} if $(x_i|y_i)_w\rightarrow \infty$ as $i\to\infty$.
\item
The {\it Gromov boundary} $\partial^* X$ of $X$ is defined to be the set of all equivalence classes of Gromov sequences, and $X^*=X \cup \partial^*  X$ is called the {\it Gromov closure} of $X$. For the description of the topology of $X^*$ we refer to \cite[Page 429]{BrHa}.
\item
For $a\in X$ and $b\in \partial^* X$, the Gromov product $(a|b)_w$ of $a$ and $b$ is defined by
$$(a|b)_w= \inf \big\{ \liminf_{i\rightarrow \infty}(a|b_i)_w\; \big|\; \{b_i\}\in b\big\}.$$
\item
For $a,$ $b\in \partial^* X$, the Gromov product $(a|b)_w$ of $a$ and $b$ is defined by
$$(a|b)_w= \inf \big\{ \liminf_{i\rightarrow \infty}(a_i|b_i)_w\;\big|\; \{a_i\}\in a\;\;{\rm and}\;\; \{b_i\}\in b\big\}.$$
\end{enumerate}
\end{enumerate}

%Note that if $(X,d)$ is proper and geodesic, its Gromov boundary is also equivalent to the geodesic boundary $\partial_G X$ (cf. \cite{BrHa}). Here, $\partial_G X$ is defined as the set of all equivalence classes of geodesic rays, where two geodesic rays are said to be {\it equivalent} if they have finite Hausdorff distance (cf. \cite{BS,BrHa}). In this paper, we always consider $\partial_G X=\partial^* X$ as a set of points.

Now, we define a metric on the boundary at infinity of a Gromov hyperbolic space via the extended Gromov products, see \cite{BrHa,BuSc}.

\bdefe\label{a-1}
Let $X$ be a $\delta$-hyperbolic space with $\delta >0$ and $w\in X$ a fixed point. For $0<\varepsilon<\min\{1,\frac{1}{5\delta}\}$, define
$$\rho_{w,\varepsilon}(\xi,\zeta)=e^{-\varepsilon(\xi|\zeta)_w}$$
for all $\xi,\zeta\in \partial^* X$ with the convention $e^{-\infty}=0$.

Let
$$d_{w,\varepsilon} (\xi,\zeta):=\inf\bigg\{\sum_{i=1}^{n} \rho_{w,\varepsilon} (\xi_{i-1},\xi_i)\;\big|\;n\geq 1,\xi=\xi_0,\xi_1,\ldots,\xi_n=\zeta\in \partial^* X\bigg\}.$$
Then $(\partial^* X, d_{w,\varepsilon})$ is a metric space with
$$\rho_{w,\varepsilon}/2\leq d_{w,\varepsilon}\leq \rho_{w,\varepsilon},$$
and we call $d_{w,\varepsilon}$ the {\it visual metric} on $\partial^* X$ with respect to $w\in X$ and the parameter $\varepsilon$.
\edefe

In \cite{BHK}, Bonk, Heinonen and Koskela  introduced the concept of rough starlikeness for Gromov hyperbolic space with respect to a given base point in the space. They also proved that both bounded uniform spaces and Gromov hyperbolic domains in $\IR^n$ are roughly starlike. It turns out that this property is very useful, see for instance \cite{BB03}.

\bdefe Let $\kappa\geq 0$. Suppose that $(X,d)$ is a proper, geodesic $\delta$-hyperbolic metric space and that $w\in X$. We say that $X$ is {\it $\kappa$-roughly starlike} with respect to $w$ if for each $x\in X$, there exists a point $\xi\in\partial^* X$ and a geodesic ray $\alpha=[w,\xi)$ satisfying
$$\dist(x,\alpha)\leq \kappa.$$
\edefe

Further, V\"{a}is\"{a}l\"{a} extended their ideas and introduced the following definition in \cite{Va05a}.
\bdefe Let $\kappa\geq 0$. Suppose that $(X,d)$ is a proper, geodesic $\delta$-hyperbolic metric space and that $\xi\in\partial^* X$. We say that $X$ is {\it $\kappa$-roughly starlike} with respect to $\xi$ if for each $x\in X$, there is a point $\eta\in\partial^* X$ and a geodesic line $\gamma=[\xi,\eta]$ joining $\xi$ and $\eta$ such that $$\dist(x,\gamma)\leq \kappa.$$
\edefe

\subsection{Bonk-Heinonen-Koskela uniformization}\label{sec-z}
We now recall the following conformal deformations of proper geodesic Gromov hyperbolic spaces that were introduced by Bonk, Heinonen and Koskela, see \cite[Chapter $4$]{BHK}. We remark that this uniformization theory has many applications, see, e.g., \cite{BB03, KLM14, ZRL}.

Let $D\subsetneq \mathbb{R}^n$ be a $\delta$-hyperbolic domain, and $k$ its quasihyperbolic metric. Fix a base point $w\in D$, and consider the family of conformal deformations of $(D,k)$ by the densities
$$ \rho_\epsilon(x)=e^{-\epsilon k(x,w)}\;\;(\epsilon>0).$$
For $x$, $y\in D$, let
 \be\label{a-5} d_\epsilon(x,y)=\inf\int_{\gamma} \rho_\epsilon \; ds_k,\ee
where $ds_k$ is the arc-length element with respect to the metric $k$ and the infimum is taken over all rectifiable curves $\gamma$ in $D$ with endpoints $x$ and $y$.

Then $d_\epsilon$ are metrics on $D$, and we denote the resulting metric spaces by $D_\epsilon=(D,d_\epsilon)$. Moreover, $k_\epsilon$ is the quasihyperbolic metric of $D_\epsilon$, $\overline{D}_\epsilon$ and $\partial_\epsilon D$ denote the metric completion and metric boundary of $D$ with respect to $d_\epsilon$, respectively.

Finally, we conclude with some auxiliary results for our needed.
%The authors are indebted to the referee for his/her suggestions about the following result and its proof.

\blem\label{r-1} Let $D\subsetneq \mathbb{R}^n$ be a $\delta$-hyperbolic domain with $n\geq 2$, and $k$ its quasihyperbolic metric. There are constants $A,\kappa,M,C,\varepsilon_0$ that depend only on $\delta$ such that for a fixed $\epsilon\in (0,\varepsilon_0]$ we have:

$(a)$ $(D,k)$ is a complete, proper and geodesic metric space.

$(b)$ $(D,k)$ is $\kappa$-roughly starlike with respect to each point of $D^*$.

$(c)$ $D_\epsilon$ is $A$-uniform and bounded $($with  diameter at most $2/\epsilon$$)$.

$(d)$ The identity map $(D,k)\to (D,k_\epsilon)$ is $M$-bilipschitz.

$(e)$ For all $x,y\in D$, we have
 $$C^{-1}d_\epsilon(x,y) \leq \epsilon^{-1}e^{-\epsilon (x|y)_w} \min\{ 1,\epsilon k(x,y) \} \leq Cd_\epsilon(x,y).$$

$(f)$  There is a natural map $D^* \xrightarrow{\varphi}\overline{D}_\epsilon$ that is a bijection.

$(g)$ There is a natural $\theta$-quasim\"obius identification $(\partial^*D,\rho)\xrightarrow{\psi} \partial_\epsilon D$ with $\psi=\varphi|_{\partial^*D}$, where $\rho$ is a visual metric on $\partial^*D$ and $\theta$ is a self-homeomorphism of $[0,\infty)$ depending only on $\delta$ and the parameter of $\rho$.
\elem
\bpf $(a)$ See \cite[Proposition 2.8]{BHK}.

$(b)$ See \cite[Theorem 3.22]{Va05a} and also \cite[Theorem 3.6]{BHK}.

$(c)$  See \cite[Proposition 4.5]{BHK}.

$(d)$ See \cite[Proposition 4.37]{BHK}.

$(e)$ See \cite[Lemma $4.10$]{BHK}.

$(f)$ See \cite[Section 2.21]{Va05a} or \cite[Section 6]{HSX} for the notion of a {\it natural} map; this is a continuous  extension of the identity map
$$D^* \supset D \xrightarrow{\id} D\subset \overline{D}_\epsilon,$$
where the topology is given as described in \cite{BrHa,Va05b}. In fact, by the fact $(e)$ or the proof of \cite[Proposition 4.13]{BHK}, we have the following:

\emph{There is a bijection $D^* \xrightarrow{\varphi}\overline{D}_\epsilon$ that satisfies $\varphi|_D=\id_D$ and a sequence $\overline{x}=\{x_i\}$ in $D$ that satisfies $d_\epsilon(x_i,\xi)\to 0$ for some $\xi\in \partial_\epsilon D$ if and only if $\overline{x}$ is a Gromov sequence in $(D,k)$ and $\varphi(\hat{x})=\xi$ where $\hat{x}\in \partial^* D$ is the equivalence class of $\overline{x}$.}

$(g)$ Note that $D_\epsilon$ is induced by the density $\rho_\epsilon(x)=e^{-\epsilon k(x,w)}$. Thus it follows from \cite[Proposition 4.13]{BHK} that there is a natural $\eta$-quasisymmetric identification $(\partial^*D,\rho_{w,\varepsilon})\to \partial_\epsilon D$, where $\rho_{w,\varepsilon}$ is a visual metric on $\partial^*D$ based at $w$ with parameter $\varepsilon$. By \cite[Corollary 5.2.9]{BuSc}, we see that $\partial^* D$ endowed with any two visual metrics are quasim\"obius equivalent. Hence there is a natural $\theta$-quasim\"obius identification $(\partial^*D,\rho)\to \partial_\epsilon D$ for all visual metric $\rho$ on $\partial^*D$, as desired.

\epf

\section{Proof of Theorem \ref{thm-1}}\label{sec-3}

Here we assume that
\begin{enumerate}
  \item $D\subsetneq \mathbb{R}^n$ is a domain,
  \item $\partial D$ is a $C$-uniformly perfect set,
  \item $f \in \mathcal{T}_K(D)=\Big\{\overline{D}  \xrightarrow{\phi} \overline{D}\;|\; \phi|_D\;\mbox{is}\;  K\mbox{-QC and}\;\phi|_{\partial D}=\id_{\partial D}\Big\}$.
\end{enumerate}

First, by Theorem \ref{Thm-A}, we know that there exists a homeomorphism $\eta:[0,\infty)\to [0,\infty)$ such that $f$ is $\eta$-quasisymmetric on $\overline{D}$ with $\eta=\eta(K,n)$. We show that for all $x\in D$,
$$j_D(x,f(x))\leq M:=2\log(1+4\eta(C)).$$

For a given point $x\in D$, take a point $x_0\in \partial D$ such that
$$|x-x_0|=d(x).$$
Next, we claim that there is a point $x_1\in \partial D$ satisfying
\be\label{z-01} \frac{d(x)}{C} \leq |x_0-x_1| \leq 4d(x). \ee
The proof of the claim is divided into two cases.

For the first case, suppose $\partial D\subseteq B(x,2d(x))$. Note that the boundary of $D$ is taken in the topology of the Riemann sphere $\mathbb{R}^n\cup \{\infty\}$. Thus $D$ is bounded with
$$\diam (\partial D)=\diam (D)\geq 2d(x).$$
Then we may pick a point $x_1\in \partial D$ with
$$|x_0-x_1|\geq \frac{1}{2} \diam (D) \geq d(x).$$
Moreover, because $\partial D\subseteq B(x,2d(x))$, we have
$$|x_0-x_1|\leq 4d(x),$$
as desired.

For the remaining case, suppose $\partial D\not\subseteq B(x,2d(x))$. A direct computation shows that
$B(x_0,d(x)) \subseteq B(x,2d(x))$, which implies
$$\partial D\setminus B(x_0,d(x))\neq \emptyset.$$
Because $\partial D$ is $C$-uniformly perfect with respect to the spherical metric, by Remark \ref{a-2} we know that $\partial D$ is uniformly perfect with respect to the Euclidean metric of $\mathbb{R}^n$. So there is no loss of generality in assuming that $\partial D$ is also $C$-uniformly perfect with the same constant in $\mathbb{R}^n$. Thus we see that there exists $x_1\in \partial D$  such that
$$\frac{d(x)}{C} \leq |x_0-x_1| \leq d(x),$$
and we obtain (\ref{z-01}).

Furthermore, we are going to show that there exists a constant $C_1\geq 1$ depending on $\eta$ and $C$ such that
\be\label{z-02}  \frac{d(x)}{C_1}\leq d(f(x))\leq |f(x)-f(x_0)|\leq C_1 d(x).   \ee
Indeed, we only need to find a constant $C_1\geq 1$ satisfying
$$d(f(x))\leq |f(x)-f(x_0)|\leq C_1 d(x),$$
because the other direction follows from a symmetric argument along with the fact that the inverse map $f^{-1}$ is $\eta'$-quasisymmetric with identity boundary values, where $\eta'(t)=\eta^{-1}(t^{-1})^{-1}$ for all $t>0$ \cite{TV}.

Because $f$ is $\eta$-quasisymmetric on $\overline{D}$ and $f|_{\partial D}=\id_{\partial D}$, we may compute from (\ref{z-01}) that
\begin{eqnarray*} |f(x)-f(x_0)| &\leq&\eta\Big(\frac{|x-x_0|}{|x_1-x_0|}\Big)|f(x_1)-f(x_0)|
\\ \nonumber&\leq& \eta(C)|x_1-x_0|
\\ \nonumber&\leq& 4\eta(C) d(x),\end{eqnarray*}
which yields (\ref{z-02}) by taking $C_1=4\eta(C)$.

Since $f|_{\partial D}=\id_{\partial D}$, by (\ref{z-02}) we find that

\begin{eqnarray*} |f(x)-x| &\leq& |f(x)-f(x_0)|+|x_0-x|
\\ \nonumber&\leq& C_1 d(x)+d(x)
\\ \nonumber&\leq& C_1(C_1+1)\min\{d(x),d(f(x))\},
\end{eqnarray*}
which implies
$$j_D(x,f(x))=\log\Big(1+ \frac{|f(x)-x|}{\min\{d(x),d(f(x))\}}\Big)\leq 2\log(1+C_1)=M,$$
as desired. \qed

\section{Proof of Theorem \ref{thm-2}}\label{sec-4}
Throughout this section, we assume that:
\begin{enumerate}
  \item $D\subsetneq \mathbb{R}^n$ is a $\delta$-hyperbolic domain,
  \item $\partial^* D$ equipped with a visual metric $\rho$ is $C$-uniformly perfect,
  \item $f \in \mathcal{T}_K^*(D)=\Big\{D \xrightarrow{\phi} D\;\big |\;\phi\;\mbox{ is }K\mbox{-QC so that}\;\partial \phi=\id_{\partial^* D} \Big\}$.
\end{enumerate}
Our goal is to find a constant $L$ such that, for all $x\in D$
$$k_D(x,f(x))\leq L.$$

Fix a point $w\in D$. Because $D$ is $\delta$-hyperbolic, we know from Lemma \ref{r-1}$(a)$ and $(b)$ that $(D,k)$ is a proper geodesic metric space, and there is a constant $\kappa=\kappa(\delta)\geq 0$ such that $D$ is $\kappa$-roughly starlike with respect to $w$.

Let $D_\epsilon=(D,d_\epsilon)$ be the BHK-uniformization of $(D,k)$ obtained via the dampening conformal deformation as described in Subsection \ref{sec-z} with $d_\epsilon$ defined as in (\ref{a-5}). We record some required auxiliary results. By Lemma \ref{r-1}, there are constants $A,M,C_1$ that depend only on $\delta$ such that

$(a)$  $D_\epsilon$ is $A$-uniform and bounded $($with  diameter at most $2/\epsilon$$)$.

$(b)$ The identity map $\vartheta:(D,k)\to (D,k_\epsilon)$ is $M$-bi-Lipschitz, where $k_\epsilon$ is the quasihyperbolic metric of $D_\epsilon$.

$(c)$ For all $x,y\in D$, we have
 $$C_1^{-1}d_\epsilon(x,y) \leq \epsilon^{-1}e^{-\epsilon (x|y)_w} \min\{ 1,\epsilon k(x,y) \} \leq C_1 d_\epsilon(x,y).$$

$(d)$  There is a natural map $D^* \xrightarrow{\varphi}\overline{D}_\epsilon$ that is a bijection.

$(e)$ There is a natural $\theta$-quasim\"obius identification $(\partial^*D,\rho)\xrightarrow{\psi} \partial_\epsilon D$ with $\psi=\varphi|_{\partial^*D}$, where $\theta$ is a self-homeomorphism of $[0,\infty)$ depending only on $\delta$ and the parameter of $\rho$.

Moreover, by the assumption that $(\partial^*D, \rho)$ is $C$-uniformly perfect, it follows from the fact $(e)$ and Lemma \ref{Thm-z0} that

$(f)$ $\partial_\epsilon D$ is  $C_0$-uniformly perfect with $C_0$ depending only on $\theta$ and $C$.

Since $f\in \mathcal{T}_K^*(D)$, we define a map $f^*:D^*\to D^*$
\beqq
 f^*:=\begin{cases}
 \ds f & \mbox{ in $D$},\\
\partial f & \;\mbox{in $\partial^*D$.}
 \end{cases}
\eeqq
Further, we consider the map
$$g:=\varphi \circ f^* \circ \varphi^{-1}:\overline{D_\epsilon} \to \overline{D_\epsilon}$$
induced by $f^*$.   For later use we need to prove the following lemma.

\blem\label{s-1} The map $g:\overline{D_\epsilon}\to \overline{D_\epsilon}$ satisfies the following properties:
\begin{enumerate}
  \item\label{zz1} $g$ is a homeomorphism with $g|_{\partial_\epsilon D}=\id_{\partial_\epsilon D}$;
  \item\label{zz3} $g|_{(D,d_\epsilon)}$ is $\eta_0$-quasisymmetric with $\eta_0=\eta_0(\delta,K,n)$.
\end{enumerate}
\elem

The proof of Lemma \ref{s-1} is divided into two parts. 

\subsection{ Proof of Lemma \ref{s-1}$(1)$.} The assertion follows from the following statements:

$(i)$ $g$ is a bijection,

$(ii)$ $g|_{D_\epsilon}$ is a homeomorphism,

$(iii)$ The continuous extension of $g|_{D_\epsilon}$ to the boundary  $\partial_\epsilon D$ satisfies $g|_{\partial_\epsilon D}=\id_{\partial_\epsilon D}$.

The statement $(i)$ follows from the fact that both $f^*$ and $\varphi$ are bijections. For $(ii)$, the statement $(d)$ implies that $\varphi$ is a natural map with $\varphi|_D=\id_D$. Thus by the definitions of $g$ and $f^*$, we know that $g|_{D_\epsilon}=f|_D$. Now it follows from \cite[Proposition 2.8]{BHK} that the identity maps $(D,|\cdot|)\to (D,k)$ and $(D,d_\epsilon)\to (D,k_\epsilon)$ are both homeomorphisms. Since $f|_{(D,|\cdot|)}$ is a homeomorphism, by the statement $(b)$, $g|_{D_\epsilon}$ is a homeomorphism as well.

It remains to show $(iii)$. For any sequence $\{x_i\}$ in $D_\epsilon$ with $d_\epsilon(x_i,p)\to 0$ for some $p\in\partial_\epsilon D$, the statement $(d)$ implies that $\{\varphi^{-1}(x_i)\}=\{x_i\}$ is a Gromov sequence in $(D,k)$ with $\varphi^{-1}(p)=\xi\in \partial^*D$.

Since $f\in \mathcal{T}_K^*(D)$, we see from Theorem \ref{Thm-1} that $f:(D,k)\to (D,k)$ is a rough quasi-isometry. Then by  \cite[Proposition 6.3]{BS}, $\{f\circ \varphi^{-1}(x_i)\}=\{f(x_i)\}$ is also Gromov in $(D,k)$ with $f^*(\xi)=\xi\in \partial^*D$. Again by the statement $(d)$, the sequence $\{g(x_i)=\{\varphi\circ f\circ \varphi^{-1}(x_i)\}\}$ satisfies
$$d_\epsilon(g(x_i),p)=d_\epsilon(\varphi\circ f\circ \varphi^{-1}(x_i),\varphi(\xi))=d_\epsilon(f(x_i),p)\to 0.$$
This ensures that the continuous extension of $g|_{D_\epsilon}$ to the boundary  $\partial_\epsilon D$ satisfies $g|_{\partial_\epsilon D}=\id_{\partial_\epsilon D}$, as required.
\qed

\subsection{ Proof of Lemma \ref{s-1}$(2)$.}
We begin with some preparations and divide the proof into several steps.

Let $D_\epsilon'=(D,d_\epsilon')$ be the BHK-uniformization of $(D,k)$ obtained via the dampening conformal
deformation as described in Subsection 2.5 with $d_\epsilon'$ defined as in (\ref{a-5}), but now we use the base point
$w'=f(w)$. Again by Lemma \ref{r-1}$(b)$, $(c)$ and $(d)$, we know that
\begin{enumerate}
  \item $(D,k)$ is $\kappa$-roughly starlike with respect to $w'$,
  \item $(D,d_\epsilon')$ is $A$-uniform, and
  \item the identity map $\vartheta': (D,k) \to (D,k_\epsilon')$ is $M$-bi-Lipschitz, where $k_\epsilon'$ is the quasihyperbolic metric of the space $(D,d_\epsilon')$.
\end{enumerate}
Moreover, by Lemma \ref{r-1}$(e)$, a direct computation shows that the identity map
$$\phi:(D,d_\epsilon) \to (D,d_\epsilon')$$
is $\theta_0$-quasim\"{o}bius with $\theta_0(t)=C't$ and $C'$ depending only on $\delta$. Then $g$ induces a map
$$h:=\phi\circ g: (D,d_\epsilon) \to (D,d_\epsilon').$$

\textbf{Outline of the proof of  Lemma \ref{s-1}$(2)$.} In the following, we first show that $h$ is quasisymmetric, see Lemma \ref{s-3}. Because the composition of a quasim\"{o}bius map and a quasisymmetric map is also quasim\"{o}bius (cf. \cite{Va85,WZ17}), we observe that $g=\phi^{-1}\circ h$ is quasim\"{o}bius. After that, we want to use Theorem \ref{Thm-2} to check the quasisymmetry of $g$. So we only need to verify the three-point condition (\ref{a-6}) stated in Theorem \ref{Thm-2}, see Lemma \ref{sz-1}.

\blem\label{s-3} The map $h:(D,d_\epsilon) \to (D,d_\epsilon')$ is $\theta$-quasisymmetric with $\theta$ depending only on $n,\delta$ and $K$.
\elem

\bpf First, we record some results from the previous arguments:
\begin{enumerate}
  \item $(D,d_\epsilon)$ and $(D,d_\epsilon')$ are both $A$-uniform and so they are $A$-quasiconvex,
  \item the identity maps $\vartheta: (D,k) \to (D,k_\epsilon)$ and $\vartheta': (D,k) \to (D,k_\epsilon')$ are both $M$-bilipschitz.
\end{enumerate}

Consider the identity maps
$$\tau:(D,|\cdot|)\to (D,d_\epsilon)$$
and
$$\tau'=\phi\circ \tau:(D,|\cdot|)\to (D,d_\epsilon').$$
Since $\varphi|_D=\phi=\id_D=\tau$, we have $h(x)=g(x)=f(x)$ for all $x\in D$.
Therefore, for all $x\in D$,
$$h(x)=\phi\circ g(x)=\phi\circ\varphi \circ f\circ \varphi^{-1}(x)=\tau' \circ f\circ \tau^{-1}(x).$$
Then we show that

\bcl\label{s-2} The map $h:(D,d_\epsilon) \to (D,d_\epsilon')$  is $q$-locally weakly $H$-quasisymmetric with $q$ and $H$ depending only on $n,\delta$ and $K$.
\ecl

By argument $(2)$ and \cite[Theorem 3.7]{HLL}, we see that the restrictions of $\tau^{-1}$ and $\tau'$ on each subdomain of $D$ are both $M'$-bilipschitz with respect to the quasihyperbolic metrics, where $M'$ depends only on $M$ and $A$. Moreover, it follows from \cite[Theorem 1.8]{HRWZ} that both $\tau^{-1}$ and $\tau'$ are $q_1$-locally weakly $H_1$-quasisymmetric on $D$, where $q_1$ and $H_1$ depend only on $M'$ and $A$.

Furthermore, because $f|_D$ is $K$-quasiconformal, we see from \cite[Theorem 11.14]{Hei} that $f|_D$ is $\frac{1}{2}$-locally $\eta_1$-quasisymmetric for some homeomorphism $\eta_1:[0,\infty) \to [0,\infty)$ with $\eta_1$ depending only on $K$ and $n$.

By \cite[Theorem 1.12]{HRWZ}, we know that the composition of two locally weakly quasisymmetric maps is also locally weakly quasisymmetric. This fact together with the locally weak quasisymmetry of $\tau',f$ and $\tau^{-1}$, shows that $h:(D,d_\epsilon) \to (D,d_\epsilon')$  is $q$-locally weakly $H$-quasisymmetric because
$h=\tau' \circ f\circ \tau^{-1}$, which shows Claim \ref{s-2}.

Note that $(D,d_\epsilon)$ and $(D,d_\epsilon')$ are both $A$-quasiconvex. To prove that $h$ is quasisymmetric, it follows from \cite[Theorem 6.6]{Va99} that we only need to find a constant $H_0\geq 1$ depending only on $n,\delta$ and $K$ such that $h$ is weakly $H_0$-quasisymmetric.

For each triple of distinct points $x,y,z\in D$ with $d_\epsilon (x,y)\leq d_\epsilon (x,z)$, we show that
\be\label{zz4} d'_\epsilon(h(x), h(y)) \leq H_0 d'_\epsilon(h(x),h(z)) \ee
for some constant $H_0$.

Let
$$t=\frac{\epsilon}{M} \log(1+q)<1.$$
Since $f|_D$ is $K$-quasiconformal, we know that the inverse map $f|_D^{-1}$ is $K'$-quasiconformal with $K'=K'(K,n)$ because the inverse map of a quasiconformal homeomorphism is again quasiconformal. Thus it follows from Theorem \ref{Thm-1} that there is a  homeomorphism $\psi=\psi(n,K):[0,\infty) \to [0,\infty)$ such that for all $u,v\in D$,
\be\label{q-0} \psi^{-1}(k(u,v)) \leq k(f(u),f(v)) \leq \psi(k (u,v)). \ee

We divide the proof of $(\ref{zz4})$  into three cases.

\bca Suppose $\epsilon k(x,z)< t$. \eca
In this case, by the elementary inequality \cite[(2.4)]{BHK} and the choice of $t$ we have
$$\frac{d_\epsilon (x,y)}{d_\epsilon (x)} \leq \frac{d_\epsilon (x,z)}{d_\epsilon (x)}\leq e^{k_\epsilon(x,z)}-1 < q,$$
which implies that $x,y,z\in B_\epsilon(x, q d_\epsilon (x))$, where $d_\epsilon (x)=d_\epsilon (x, \partial_\epsilon D)$ and
$$B_\epsilon(x, q d_\epsilon (x))=\{y\in D\;|\;d_\epsilon(x,y)<q d_\epsilon (x)\}.$$
Hence we obtain $(\ref{zz4})$ from Claim \ref{s-2} by choosing $H_0=H$.

\bca Suppose $\epsilon k(x,z)\geq t$ and $\epsilon k(x,y)<1$.\eca

Because $\epsilon k(x,z)\geq t$, we see from (\ref{q-0}) that
\be\label{q-1} \epsilon k(h(x),h(z)) =\epsilon k(f(x),f(z))  \geq \epsilon \psi^{-1}(t/\epsilon )=:1/C_2, \ee
and similarly
\be\label{q-2} \epsilon k(h(x),h(y))=\epsilon k(f(x),f(y))  \leq \epsilon \psi(1/\epsilon)=:C_3.\ee
Then, by using (\ref{q-1}), (\ref{q-2}) and the statement $(c)$, we have
\beq\nonumber
\frac{d'_\epsilon(h(x), h(y))}{d'_\epsilon(h(x),h(z))} & \leq & C_1^2 e^{\epsilon (h(x)|h(z))_{w'}- \epsilon (h(x)|h(y))_{w'}} \frac{\min\{1, \epsilon k(h(x),h(y))\}}{\min\{1,\epsilon k(h(x),h(z))\}}
\\ \nonumber & \leq & C_1^2 C_2 e^{\epsilon k(h(x),h(y))}
\\ \nonumber & \leq & C_1^2 C_2 e^{C_3},
\eeq
where $C_1=C_1(\delta)$ is the constant of the statement $(c)$. By setting $H_0=C_1^2 C_2 e^{C_3}$, we obtain  $(\ref{zz4})$.

\bca Suppose $\epsilon k(x,z)\geq t$ and $\epsilon k(x,y)\geq 1.$\eca
Because $d_\epsilon (x,y)\leq d_\epsilon (x,z)$, we note again from the statement $(c)$ that
$$e^{\epsilon (x|z)_w-\epsilon (x|y)_w} \frac{\min\{1, \epsilon k(x,y)\}}{\min\{1, \epsilon k(x,z)\}}\leq C_1^2$$
and so
\be\label{zz5} (x|z)_w- (x|y)_w\leq \frac{2\log C_1}{\epsilon}. \ee
It follows from  (\ref{q-0}) that $f:(D,k)\to (D,k)$ and its inverse map are both $\psi$-uniformly continuous (for the definition see \cite[Section 2]{Va99}). Because $\epsilon k(x,z)\geq t$, by (\ref{q-0}) we obtain
\be\label{zz6}\epsilon k(h(x),h(z))\geq 1/C_2.\ee

Moreover, as $(D,k)$ is geodesic, by \cite[Theorem 2.5]{Va99} we know that there are positive constants $\lambda$ and $\mu$ depending only on $\psi$ such that $f:(D,k)\to (D,k)$ is $(\lambda,\mu)$-rough quasi-isometry. Consequently, we see from $(\ref{zz5})$ and  \cite[Proposition 5.5]{BS} that there is a constant $C_4=C_4(\lambda, \mu, \delta, C_1,\epsilon)$ such that
$$\epsilon (h(x)|h(z))_{w'}- \epsilon (h(x)|h(y))_{w'}\leq C_4.$$
So again by the statement $(c)$ and (\ref{zz6}), we obtain
\beq\nonumber
\frac{d'_\epsilon(h(x), h(y))}{d'_\epsilon(h(x),h(z))} & \leq & C_1^2 e^{\epsilon (h(x)|h(z))_{w'}- \epsilon (h(x)|h(y))_{w'}} \frac{\min\{1, \epsilon k(h(x),h(y))\}}{\min\{1,\epsilon k(h(x),h(z))\}}
\\ \nonumber & \leq & C_1^2 C_2 e^{C_4},
\eeq
as needed.
\epf

As mentioned before, to prove the quasisymmetry of $g$, it suffices to show that $g$ satisfies the three-point condition (\ref{a-6}) stated in Theorem \ref{Thm-2}.

\blem\label{sz-1} There are three distinct points $\xi_1,\xi_2,\xi_3\in \partial_\epsilon D$, and a number $\lambda_0>0$ such that
$$d_\epsilon (g(\xi_i),g(\xi_j))=d_\epsilon (\xi_i,\xi_j) \geq \lambda_0 \diam_\epsilon( D_\epsilon),$$
for all $i\neq j\in \{1,2,3\}$.
\elem
\bpf
To this end, we first show that
\be\label{a-3} \diam_\epsilon (\partial_\epsilon D)\leq \diam_\epsilon (D_\epsilon) \leq M_1 \diam_\epsilon (\partial_\epsilon D),\ee
where $M_1=A e^{\epsilon\kappa}$.

Because $\diam_\epsilon (\partial_\epsilon D)\leq \diam_\epsilon (D_\epsilon)$, it suffices to check the second inequality.

We note from the statement $(d)$ that there is a natural identification $\varphi:\partial^*D\to \partial_\epsilon D.$
Fix a point $\xi'\in \partial_\epsilon D$ and take $\xi\in \partial^*D$ with $\varphi(\xi)=\xi'$. By Lemma \ref{r-1}$(b)$, $(D,k)$ is $\kappa$-roughly starlike with respect to $\xi$. It follows that for $w\in D$, there is another point $\zeta \in \partial^*D$ and a quasihyperbolic geodesic line $\gamma=[\xi,\zeta]_k$ joining $\xi$ and $\zeta$ such that $\{\gamma(-n)\}_{n=1}^\infty\in\xi$, $\{\gamma(n)\}_{n=1}^\infty\in\zeta$ and the quasihyperbolic distance
$$\dist_k(w, [\xi,\zeta]_k) \leq \kappa.$$
This shows that there is a point $w_0\in [\xi,\zeta]_k$ satisfying $k(w,w_0)\leq \kappa.$

Let $(\mathbb{R},|\cdot|) \xrightarrow{\gamma} (D,k)$ be a $k$-arclength parametrization of the quasihyperbolic geodesic line $(\xi,\zeta)_k$ with $\gamma(0)=w_0$, $\gamma(-\infty)=\xi$ and $\gamma(+\infty)=\zeta$. Then for each $x\in (\xi,\zeta)_k$,
$$k(x,w)\leq k(x,w_0)+k(w_0,w)\leq \kappa+k(x,w_0).$$
Therefore,
$$\ell_\epsilon(\gamma)=\int_{\gamma} e^{-\epsilon k(w,x)} \;ds_k\geq e^{-\epsilon\kappa}\int_{\gamma} e^{-\epsilon k(w_0,x)} \; ds_k= e^{-\epsilon\kappa}\int_{-\infty}^{+\infty} e^{-\epsilon|t|}\; dt= 2e^{-\epsilon\kappa}\epsilon^{-1}.$$

Moreover, by the proof of \cite[Proposition 4.5]{BHK}, $\gamma$ is an $A$-uniform arc in $D_\epsilon$. This gives
$$\diam_\epsilon(\partial_\epsilon D)\geq d_\epsilon(\xi,\zeta)\geq A^{-1} \ell_\epsilon(\gamma)\geq \frac{e^{-\epsilon\kappa}}{A} \frac{2}{\epsilon}\geq \frac{e^{-\epsilon\kappa}}{A}\diam_\epsilon ({D}_\epsilon),$$
where the last inequality follows from the statement $(a)$ that $\diam_\epsilon (D_\epsilon) \leq 2/\epsilon.$
Hence we obtain (\ref{a-3}).

Next, take two points $\xi_1$ and $\xi_2$ in $\partial_\epsilon D$ with
$$d_\epsilon (\xi_1,\xi_2)=\diam_\epsilon( \partial_\epsilon D).$$
The statement $(f)$ guarantees that $\partial_\epsilon D$ is $C_0$-uniformly perfect. Because $\partial_\epsilon D\setminus B_\epsilon (\xi_1, \frac{1}{2}d_\epsilon (\xi_1,\xi_2))\neq \emptyset$, it follows that there is a point $\xi_3\in \partial_\epsilon D$ with
$$\frac{d_\epsilon (\xi_1,\xi_2)}{2C_0} \leq d_\epsilon (\xi_1,\xi_3) \leq \frac{1}{2} d_\epsilon (\xi_1,\xi_2)$$
and therefore
$$d_\epsilon (\xi_3,\xi_2)\geq d_\epsilon (\xi_1,\xi_2)-d_\epsilon (\xi_1,\xi_3)\geq \frac{1}{2}d_\epsilon (\xi_1,\xi_2).$$
This, together with (\ref{a-3}), shows that for all $i\neq j\in \{1,2,3\}$,
$$d_\epsilon (\xi_i,\xi_j) \geq  \frac{\diam_\epsilon( \partial_\epsilon D)}{2C_0}\geq \frac{\diam_\epsilon( D_\epsilon)}{2C_0Ae^{\epsilon\kappa}}.$$
 Hence $\lambda_0=(2C_0Ae^{\epsilon\kappa})^{-1}$ is the needed.
\epf

In the following, we continue the proof of Theorem \ref{thm-2}.

Note that we only need to check that there is a constant $\Lambda\geq 0$ depending only on $n,K,\delta$ and $C$ such that for all $x\in D$,
\be\label{t-1} k_\epsilon(x,g(x))\leq \Lambda.\ee

Indeed, we see from the statement $(b)$ and (\ref{t-1}) that
$$k(x,f(x))\leq M k_\epsilon(x,g(x))\leq \Lambda M,$$
which is the required estimate in Theorem \ref{thm-2} with the choice of $L=M \Lambda$. Thus, it remains to prove (\ref{t-1}).

\subsection{ Proof of (\ref{t-1}).}

The following arguments for (\ref{t-1}) are similar to the proof of Theorem \ref{thm-1}. For completeness, we show the details. Before proceeding further, we need to prove some technical statements.

Fix $x\in D$ and take $x_0\in \partial_\epsilon D$ such that
$$d_\epsilon(x)=d_\epsilon (x,\partial_\epsilon D)=d_\epsilon(x,x_0).$$
We first show that there is a constant $M_2=M_2(M_1,C_0)\geq 1$ and a point $x_1\in \partial_\epsilon D$ satisfying
\be\label{z1} \frac{1}{M_2}d_\epsilon(x)\leq d_\epsilon(x_1,x_0) \leq 4d_\epsilon(x).\ee

We consider two possibilities. If $\partial_\epsilon D\subseteq B_\epsilon(x,2d_\epsilon(x))$, then there is a point $x_1\in \partial_\epsilon D$ with
$$d_\epsilon(x_1,x_0) \geq \frac{1}{2} \diam_\epsilon( \partial_\epsilon D) \geq \frac{1}{2M_1} \diam_\epsilon (D_\epsilon)\geq \frac{d_\epsilon(x)}{2M_1},$$
where the penultimate inequality follows from (\ref{a-3}). Moreover, we have
$$d_\epsilon(x_1,x_0) \leq d_\epsilon(x_1,x) +d_\epsilon(x,x_0) \leq 4d_\epsilon(x),$$
which implies (\ref{z1}).

If $\partial_\epsilon D \not\subseteq B_\epsilon(x,2d_\epsilon(x))$, then we have
$$\partial_\epsilon D\setminus B_\epsilon(x_0,d_\epsilon(x)) \neq \emptyset,$$
because $B_\epsilon(x_0,d_\epsilon(x)) \subseteq B_\epsilon(x,2d_\epsilon(x))$. Moreover, as $\partial_\epsilon D$ is $C_0$-uniformly perfect, it follows that there exists some point $x_1\in \partial_\epsilon D$ such that
$$\frac{d_\epsilon(x)}{C_0}\leq d_\epsilon(x_1,x_0) \leq d_\epsilon(x),$$
as desired. Hence we obtain (\ref{z1}) with the choice of $M_2=2M_1C_0$.

Next, we show that there is a constant $M_3\geq 1$ such that
\be\label{z2} \frac{d_\epsilon(x)}{M_3}\leq d_\epsilon (g(x)) \leq d_\epsilon(g(x),g(x_0)) \leq M_3 d_\epsilon(x). \ee

By symmetry, we only need to show that
$$d_\epsilon(g(x),g(x_0)) \leq M_3 d_\epsilon(x),$$
because the inverse map of an $\eta$-quasisymmetric homeomorphism is $\eta'$-quasisymmetric with  $\eta'(t)=\eta^{-1}(t^{-1})^{-1}$ for all $t>0$ (cf. \cite{TV}).

It follows from Lemma \ref{s-1} and \cite[Theorem 2.25]{TV} that $g: \overline{D_\epsilon} \to \overline{D_\epsilon}$ is $\eta_0$-quasisymmetric with $g|_{\partial_\epsilon D}=\id_{\partial_\epsilon D}$, where $\eta_0=\eta_0(\delta,K,n,C)$. Now by (\ref{z1}), we have
\beq\nonumber
d_\epsilon(g(x),g(x_0))  & \leq & \eta_0\Big(\frac{d_\epsilon (x,x_0)}{d_\epsilon (x_1,x_0)}\Big) d_\epsilon(g(x_1),g(x_0))
\\ \nonumber & \leq & \eta_0(M_2) d_\epsilon(x_1,x_0)
\\ \nonumber & \leq & 4\eta_0(M_2) d_\epsilon(x),
\eeq
which shows (\ref{z2}) by taking $M_3= 4\eta_0(M_2)$.

Because $g|_{\partial_\epsilon D}=\id_{\partial_\epsilon D}$, we obtain from (\ref{z2}) that
\beq\label{zz7}
d_\epsilon(g(x),x)  & \leq & d_\epsilon(x_0,x) + d_\epsilon(g(x),g(x_0))
\\ \nonumber & \leq&  d_\epsilon(x)+M_3 d_\epsilon(x)
\\ \nonumber & \leq & M_3(M_3+1) \min\{ d_\epsilon(x), d_\epsilon(g(x))\}.
\eeq
Moreover, as $(D,d_\epsilon)$ is $A$-uniform, we see from Lemma \ref{Thm-z1} and (\ref{zz7}) that
\beq\nonumber
k_\epsilon(x,g(x))  & \leq & 4A^2\log \Big(1+\frac{d_\epsilon(g(x),x)}{\min\{ d_\epsilon(x), d_\epsilon(g(x))\}}\Big)
\\ \nonumber & \leq&  4A^2 \log [1+M_3(M_3+1)]=:\Lambda,
\eeq
as desired. This proves (\ref{t-1}).
\qed

\section{Proofs of Corollaries \ref{cor-1}, \ref{cor-2} and \ref{cor-3}}\label{sec-5}

\subsection{Proof of Corollary \ref{cor-1}} Assume that $D\subsetneq \mathbb{R}^n$ is a $\psi$-uniform domain, that $\partial D$ is a $C$-uniformly perfect set, and that  $f \in \mathcal{T}_K(D)$. We first see from Theorem \ref{thm-1} that there is a constant $M=M(n,C,K)$ such that for all $x\in D$,
$$j_D(x,f(x))\leq M.$$
This yields
\beq\nonumber
k_D(x,f(x))  & \leq & \psi(r_D(x,y))
\\ \nonumber & =&  \psi[e^{j_D(x,f(x))}-1]
\\ \nonumber &\leq&\psi(e^{M}-1)=:M',
\eeq
because $D$ is $\psi$-uniform.
\qed

\subsection{Proof of Corollary \ref{cor-2}} Assume that $D\subsetneq \mathbb{R}^n$ is a $\delta$-hyperbolic domain, that $\partial^* D$ is a $C$-uniformly perfect set, and that $f \in \mathcal{T}_K^*(D)$. By Theorem \ref{thm-2} we know that there is a constant $L=L(n,\delta,C,K)$ such that for all $x\in D$,
$$k_D(x,f(x))\leq L.$$
Then for all $x,y\in D$, we obtain
$$|k_D(f(x),f(y))-k_D(x,y)|  \leq  k_D(x,f(x))+k_D(y,f(y))\leq 2L=:L',$$
as desired. \qed

\subsection{Proof of Corollary \ref{cor-3}} Assume that $D\subsetneq \mathbb{R}^n$ is an inner $A$-uniform domain, that $\partial_I D$ is a $C$-uniformly perfect set, and that $f\in \mathcal{T}_{K}(D_I)$.

Note first that $(D,d_I)$ is an $A$-uniform incomplete locally compact metric space. It follows from \cite[Theorem 3.6]{BHK} that $(D,k)$ is a proper and geodesic $\delta$-hyperbolic space with $\delta=\delta(A)\geq 0$, i.e., $D$ is a $\delta$-hyperbolic domain in $\mathbb{R}^n$.

Next, we see from \cite[Theorem 6.2]{HSX} that there is a natural map $\phi:D^*\to \overline{D}_I$
such that $\phi|_D=\id_D$ and $\phi:(\partial^*D,d_{w,\varepsilon})\to (\partial_I D,d_I)$ is  $\eta$-quasim\"obius with $\eta=\eta(A)$, where $\partial^*D$ is the Gromov boundary of $(D,k)$ and $d_{w,\varepsilon}$ is a visual metric on $\partial^*D$ with base point $w\in D$ and parameter $\varepsilon=\varepsilon(A)>0$.   Then by Lemma \ref{Thm-z0}, we see that $\partial^*D$ is $C_0$-uniformly perfect with constant $C_0=C_0(C,\eta)=C_0(C,A)$.

Now we define
$$g:=\phi\circ f\circ \phi^{-1}:D^*\to D^*.$$
As $f\in \mathcal{T}_{K}(D_I)$ and $\phi|_D=\id_D$, we know that $g\in \mathcal{T}^*_K(D).$ Consequently, by Theorem \ref{thm-2}, we see that there is a constant $H=H(n,\delta,C_0,K)=H(n,A,C,K)$ such that for all $x\in D$,
$$k_D(x,f(x))=k_D(x,g(x))\leq H,$$
as desired.
\qed

%%%%%%%%%%%%%%%%%%%
%%%%%%%%%%%%%%%%%%%
\bigskip

{\bf Acknowledgement.} The authors would like to thank Professors M. Vuorinen and Y. Li for many valuable comments and suggestions. Authors are indebted to the anonymous referee for his/her remarks that helped to improve this manuscript.

%%%%%%%%%%%%%%%%
%%%%%%%%%%%%%%%%%%

\end{document}